\documentclass[a4paper,twoside,12pt]{article}
\usepackage{amssymb} 
\begin{document}
\title{Vertex Splitting and Upper Embeddable Graphs  \footnotemark[2]
\author{Guanghua Dong$^{1,2*}$, Ning Wang$^{3}$, Yuanqiu Huang$^{1}$, Han Ren$^{4}$, Yanpei Liu$^{5}$\\
{\small\em 1.Department of Mathematics, Normal
University of Hunan, Changsha, 410081, China}\\
\hspace{-1mm}{\small\em 2.Department of Mathematics, Tianjin
Polytechnic
University, Tianjin, 300160, China}\\
\hspace{-5mm}{\small\em 3.Department of Information Science and
Technology,
Tianjin University of Finance }\\
\hspace{-74mm} {\small\em and Economics, Tianjin, 300222, China}\\
\hspace{-1mm}{\small\em 4.Department of Mathematics, East China
Normal University, Shanghai, 200062,China}\\
\hspace{-7mm}{\small\em 5.Department of Mathematics, Beijing
Jiaotong
University, Beijing, 100044, China}\\
}}\footnotetext[2]{\footnotesize \em   This work was partially
Supported  by the New Century Excellent Talents in University (Grant
No: NCET-07-0276 (Y. Huang)), the National Natural Science
Foundation of China (Grant No. 11171114 (H. Ren); 10871021 (Y.
Liu)), and the China Postdoctoral Science Foundation funded project
(Grant No: 20110491248 (G. Dong)). }
 \footnotetext[1]{\footnotesize \em E-mail: gh.dong@163.com(G. Dong); ninglw@163.com(N. Wang); hyqq@hunnu.edu.cn(Y. Huang); hren@math.ecnu.edu.cn(H. Ren);
ypliu@bjtu.edu.cn(Y. Liu). }

\date{}
\maketitle

\vspace{-.8cm}
\begin{abstract}
The $weak$ $minor$ $\underline{G}$ of a graph $G$ is the graph
obtained from $G$ by a sequence of edge-contraction operations on
$G$. A $weak$-$minor$-$closed$ family of upper embeddable graphs is
a set $\mathcal {G}$ of upper embeddable graphs that for each graph
$G$ in $\mathcal {G}$, every weak minor of $G$ is also in $\mathcal
{G}$. Up to now, there are few results providing the necessary and
sufficient conditions for characterizing upper embeddability of
graphs. In this paper, we studied the relation between the vertex
splitting operation  and the upper embeddability of graphs; provided
not only a necessary and sufficient condition for characterizing
upper embeddability of graphs, but also a way to construct
weak-minor-closed family of upper embeddable graphs from the bouquet
of circles; extended a result in $J.$ $Graph$ $Theory$ obtained by
L. Nebesk\'{y}. In addition, the algorithm complex of determining
the upper embeddability  of a graph can be reduced much by the
results obtained in this paper.


\hspace{1.77cm} flexible-edge 
\end{abstract}


\bigskip
\noindent {\bf 1. Introduction}

\bigskip

Graphs considered here are all connected, undirected, and with
minimum degree at least three. In addition, multiple edges and loops
are permitted. Terminologies and notations not defined here can be
seen in \cite{west}. The reader is assumed to be familiar with
topological graph theory, which can be find more details in
\cite{moh}, \cite{gro} or \cite{liu1}.

A graph is denoted by $G$ = ($V(G), E(G)$), and  $V(G)$, $E(G)$
denotes its vertex set and   edge set respectively. The number
$|E(G)|$ $-$ $|V(G)|$ + 1 is known as the \emph{Betti number} (or
\emph{cycle rank}) of the connected graph \emph{G},  and is denoted
by $\beta( \emph{G} )$. A $u,v$-$path$ is a path whose vertices of
degree 1 (its endpoints) are $u$ and $v$. Let \emph{T} be a spanning
tree of a connected graph \emph{G}. Define the \emph{deficiency }
$\xi( G, T )$ of a spanning tree $T$ in a graph \emph{G} to be the
number of components of $G-E(T)$ which have odd size. The deficiency
$\xi(G)$ of a graph \emph{G} is defined to be the minimum value of
$\xi( G,T )$ over all spanning tree \emph{T} of \emph{G}, $i.e.$,
$\xi(G) = min \{ \xi( G,T )\mid$ T is an spanning tree of G\}. A
$splitting$ $tree$ of a connected graph $G$ is a spanning tree $T$
for $G$ such that at most one component of $G-E(T)$ has odd size.
Let $v$ be a vertex of $G$, and $N_{G}(v)$ be the set of vertices in
$G$ adjacent to $v$, then the subgraph induced by $N_{G}(v)$ is
referred to as the $v$-$local$ subgraph, and is denoted by
$G_{loc}(v)$. The $vertex$ $splitting$ on a vertex $v$, whose degree
deg$_{G}(v)\geqslant 4$, is the replacement of the vertex $v$ by
adjacent vertices $v'$ and $v''$ and the replacement of each edge
$e=vu$ incident to $v$ either by the edge $v'u$ or by the edge
$v''u$, and the edge $v'v''$ in the new $G^{*}$ is called the
$splitting$-$edge$. If $G^{*}$ is a graph obtained from $G$ by a
vertex splitting operation on the vertex $v \in V(G)$, then the
subgraph of $G^{*}$, which is induced by $v'$, $v''$ and the
vertices adjacent to $v'$ and $v''$, is refereed to as the
$v$-$spliting$ $subgraph$ and is denoted by $G^{*}_{spl}(v)$. The
$intersection$ of two graphs $G_1$ and $G_2$ is defined as $G_1 \cap
G_2 = (V(G_1) \cap V(G_2), E(G_1) \cap E(G_2))$, and the $union$ of
$G_1$ and $G_2$ is defined as $G_1 \cup G_2 = (V(G_1) \cup V(G_2),
E(G_1) \cup E(G_2))$. A $partial$ $order$ $\mathcal {R}$ on a set
$X$ is a binary relation that is reflexive, antisymmetric, and
transitive. A $poset$, which is  short for $partially$ $ordered$
$set$, is a pair ($X; \mathcal {R}$) where $X$ is a set and
$\mathcal {R}$ is a $partial$ $order$ $relation$ on $X$. The $weak$
$minor$ $\underline{G}$ of a graph $G$, which is denoted by
$\underline{G}\preccurlyeq G$, is the graph obtained from $G$ by a
sequence of edge-contraction operations on $G$. Furthermore, a graph
$G$ is a weak minor of itself. For example, both $G_{1}$ in Fig.2
and $G_{2}$ in Fig.3 are a weak-minor of the graph $G$ in Fig.1. A
$weak$-$minor$-$closed$ family of upper embeddable graphs is a set
$\mathcal {G}$ of upper embeddable graphs that for each graph $G$ in
$\mathcal {G}$, every weak minor of $G$ is also in $\mathcal {G}$.
Obviously, the binary relation $weak$ $minor$, which is denoted by
$\preccurlyeq$, is a $partial$ $order$.

\begin{footnotesize}

\setlength{\unitlength}{0.8mm}
\begin{center}
\begin{picture}(100,20)

\put(-20,-19){\begin{picture}(10,10)

\put(22,20){\circle*{1.5}}    

\put(2,20){\circle*{1.5}}    

\put(12,40){\circle*{1.5}} 

\put(12,28.5){\circle*{1.5}} 

\put(2.5,20){\line(1,0){19}}  

\put(2.3,20.5){\line(5,4){9.2}}  

\put(2.1,20.8){\line(1,2){9.3}} 

\put(12,29){\line(0,1){10.5}}   

\put(12.8,28){\line(6,-5){9.2}}    

\put(22,20.6){\line(-1,2){9.4}}    

\put(4,14){\scriptsize {\bf {Fig.1}}:  G}

\end{picture}}

\put(40,-19){\begin{picture}(10,10)

\put(22,20){\circle*{1.5}}    

\put(2,20){\circle*{1.5}}    

\put(12,40) {\circle*{1.5}} 

\put(2.5,20){\line(1,0){19}}  

\put(2.1,20.8){\line(1,2){9.3}} 

\put(22,20.6) {\line(-1,2){9.4}}  

\qbezier(11.9,39.5)(11,24)(2.6,20.6)

\qbezier(12.1,39.5)(13,24)(21.5,20.6)

\put(4,14){{\scriptsize {\bf {Fig.2}}: $G_{1}$}}
\end{picture}}


\put(100,-19){\begin{picture}(10,10)

\put(12.1,28.3){\circle*{1.8}} 

\qbezier(2,20)(2,25)(11.3,28.3)  

\qbezier(2,20)(6,17)(11.8,27.8)  

\qbezier(12,40)(8,39)(11.5,28.9)  

\qbezier(12,40)(17,39)(12.7,29)  

\qbezier(22,20)(24,25)(12.5,28.5)  

\qbezier(22,20)(16.9,17.9)(12,27.8)  

 \put(4.5,14){{\scriptsize {\bf {Fig.3}}: $G_{2}$}}

\end{picture}}

\end{picture}
\end{center}
\end{footnotesize}

The \emph{maximum genus} $\gamma_M(G)$ of a connected graph \emph{G}
is the maximum integer \emph{k} such that there exists an embedding
of $G$ into the orientable surface of genus $k$. A graph $G$ is said
to be \emph{upper embeddable} if $\gamma_M(\emph{G})$ =
$\lfloor\frac{\beta(G)}{2}\rfloor$. Nordhaus, Stewart and White
\cite{nor} introduced the idea of the maximum genus of graphs in
1971. From then on, many interesting results have being made, mainly
concerned with the relation between the maximum genus and other
graph parameters as diameter, face size, connectivity, girth, etc.,
and the readers can find more details in
\cite{tho}\cite{chen}\cite{hua}\cite{ouy}\cite{ren}\cite{cheny}\cite{lid}\cite{neb}\cite{sko}\cite{xuo}
etc.. But few papers have provided the informations about the
problems as: (I) the relation between the upper embeddability and
vertex splitting; (II) the weak-minor-closed family of upper
embeddable graphs. The following is the details for the two
problems.

Problem I: \ Let $G$ be an upper embeddable graph,  $v$ be a vertex
of $G$ with degree no less than 4, and $G^{*}$ be the graph obtained
from $G$ through a vertex splitting operation on $v$, then $G^{*}$
may be upper embeddable or not. For example, both the graph $G_{1}$
in Fig.5 and the graph $G_{2}$  in Fig.6 are obtained from an upper
embeddable $G$ in Fig.4 through a vertex splitting operation on $v$
in $G$. The graph $G_{1}$ is upper embeddable, but $G_{2}$ is not
upper embeddable. So, a question is naturally raised: How does an
upper embeddable graph remain the upper embeddability after the
vertex splitting operation on some vertex $v$ of this graph?

\begin{footnotesize}

\setlength{\unitlength}{0.6mm}
\begin{center}
\begin{picture}(72,27.5)

\put(-68,-14){\begin{picture}(1,1)


\put(10,40){\circle*{1.5}}

\put(10,33){\circle*{1.5}}

\put(10,26){\circle*{1.5}}

\put(10,19){\circle*{1.5}}


\put(50,40){\circle*{1.5}}

\put(50,33){\circle*{1.5}}

\put(50,26){\circle*{1.5}}

\put(50,19){\circle*{1.5}}


\put(17,36.3){\circle*{1.5}}

\put(24,32.6){\circle*{1.5}}

          \put(30,29){\circle{2}}  
          \put(27.5,31.5){ $v$ }

\put(24,25.5){\circle*{1.5}}

\put(17,22){\circle*{1.5}}


\put(43,36.3){\circle*{1.5}}

\put(36,32.6){\circle*{1.5}}

\put(36,26.1){\circle*{1.5}}

\put(44,22){\circle*{1.5}}


\put(10,40){\line(0,-1){21}}

\put(10.2,39.5){\line(2,-1){19}}

 \put(10.2,18.5){\line(2,1){19}}


\put(24,25.5){\line(0,1){8}}

\put(17,22){\line(-2,1){7}}

 \put(10,33){\line(2,1){6}}


\put(10.7,40){\line(1,0){39}}

\put(50,40){\line(0,-1){20}}

 \put(10,18.4){\line(1,0){40}}


\put(50.5,39.9){\line(-2,-1){20}}

\put(50,19){\line(-2,1){19}}


\put(36,32.6){\line(0,-1){7}}

\put(44,22){\line(5,3){6}}

 \put(43,36.3){\line(2,-1){6}}

\put(22,7){{\scriptsize {\bf {Fig.4}}: G} }

\end{picture}
}



\put(75,-14){\begin{picture}(1,1)


\put(10,40){\circle*{1.5}}

\put(10,33){\circle*{1.5}}

\put(10,26){\circle*{1.5}}

\put(10,19){\circle*{1.5}}


\put(60,40){\circle*{1.5}}

\put(60,33){\circle*{1.5}}

\put(60,26){\circle*{1.5}}

\put(60,19){\circle*{1.5}}


\put(17,36.3){\circle*{1.5}}

\put(24,32.6){\circle*{1.5}}

          \put(30,29){\circle*{1.5}}  

          \put(27.5,31){ $v'$ }

\put(24,25.5){\circle*{1.5}}

\put(17,22){\circle*{1.5}}


\put(53,36.3){\circle*{1.5}}

\put(46,32.6){\circle*{1.5}}

          \put(40,29){\circle*{1.5}} 

          \put(36,31){ $v''$ }

\put(46,26.1){\circle*{1.5}}

\put(54,22){\circle*{1.5}}


\put(10,40){\line(0,-1){21}}

\put(10.2,39.5){\line(2,-1){20}}

 \put(10.2,18.5){\line(2,1){20}}


\put(24,25.5){\line(0,1){8}}

\put(17,22){\line(-2,1){7}}

 \put(10,33){\line(2,1){6}}


\put(10.7,40){\line(1,0){49}}

\put(60,40){\line(0,-1){20}}

 \put(10,18.4){\line(1,0){50}}

\put(30,29){\line(1,0){9}}

\put(60.5,39.9){\line(-2,-1){21}}

\put(60,19){\line(-2,1){20}}


\put(46,32.6){\line(0,-1){7}}

\put(54,22){\line(5,3){6}}

 \put(53,36.3){\line(2,-1){6}}

\put(26,7){{\scriptsize {\bf {Fig.6}}: $G_2$ }}
\end{picture}
}


\put(5,-14){\begin{picture}(1,1)

\put(10,40){\circle*{1.5}}

\put(10,33){\circle*{1.5}}

\put(10,26){\circle*{1.5}}

\put(10,19){\circle*{1.5}}


\put(50,40){\circle*{1.5}}

\put(50,33){\circle*{1.5}}

\put(50,26){\circle*{1.5}}

\put(50,19){\circle*{1.5}}


\put(17,36.3){\circle*{1.5}}

\put(24,32.6){\circle*{1.5}}

          \put(30,32.6){\circle*{1.5}}  

          \put(27,34.3){ $v'$ }

          \put(30,25.5){\circle*{1.5}}  

          \put(26,20.3){ $v''$ }

\put(24,25.5){\circle*{1.5}}

\put(17,22){\circle*{1.5}}


\put(43,36.3){\circle*{1.5}}

\put(36,32.6){\circle*{1.5}}

\put(36,25.5){\circle*{1.5}}

\put(44,22){\circle*{1.5}}


\put(10,40){\line(0,-1){21}}

\put(10.2,39.5){\line(2,-1){14}}

 \put(10.2,18.5){\line(2,1){14}}


\put(24,25.5){\line(0,1){8}}

\put(17,22){\line(-2,1){7}}

 \put(10,33){\line(2,1){6}}


\put(10.7,40){\line(1,0){39}}

\put(50,40){\line(0,-1){20}}

 \put(10,18.4){\line(1,0){40}}


\put(50.5,39.9){\line(-2,-1){14}}

\put(50,19){\line(-2,1){14}}


\put(36,32.6){\line(0,-1){7}}

\put(44,22){\line(5,3){6}}

 \put(43,36.3){\line(2,-1){6}}

\put(30,32.6){\line(0,-1){7}}

\put(30,32.6){\line(1,0){6}}

\put(30,32.6){\line(-1,0){6}}

\put(30,25.5){\line(-1,0){6}}

\put(30,25.5){\line(1,0){7}}


\put(22,7){{\scriptsize {\bf {Fig.5}}: $G_{1}$ }}
\end{picture}
}

\end{picture}
\end{center}
\end{footnotesize}

Problem II: \ In general, a class of upper embeddable graphs is not
closed under minors. For example, although the graph $G$ depicted in
Fig.8 is upper embeddable, the graph $G_{1}$ in Fig.7, which is a
minor of $G$, is not upper embeddable. But, if $G$ is an upper
embeddable graph then every weak minor $\underline{G}$ of $G$ is
also upper embeddable. So we can easily get a poset $\mathcal {F}$,
which is a weak-minor closed family of upper embeddable graphs, from
$G$ through a sequence of edge-contraction operations on $G$.
Obviously, the bouquet of circles $B_{\beta(G)}$, which consists of
a single vertex with $\beta(G)$ loops incident to this vertex, is
the smallest element of $\mathcal {F}$, $i.e.$, every upper
embeddable graph with $\beta(G)$ co-tree edges has bouquet circles
$B_{\beta(G)}$ as its weak-minor. However, from the example in
Fig.4-Fig.6 we can get that the bouquet circles $B_{\beta(G)}$ may
also be a weak-minor of a graph $G$ which is not upper embeddable.
So, how to get a poset $\mathcal {F}$, which is a weak-minor-closed
family of upper embeddable graphs, from the bouquet of circles
$B_{n}$ or other upper embeddable graph via series of
vertex-splitting operations on it is the second problem.

\begin{footnotesize}
\setlength{\unitlength}{0.6mm}
\begin{center}
\begin{picture}(70,13)

\put(-55,-11){\begin{picture}(30,40)

\put(18,26){\circle*{1.5}}

\put(18,14){\circle*{1.5}}

\put(30,20){\circle*{1.5}}

\put(45,20){\circle*{1.5}}

\put(57,26){\circle*{1.5}}

\put(57,14){\circle*{1.5}}

\put(30,20){\line(1,0){15}}

\put(30,20){\line(-2,1){12}}

\put(30,20){\line(-2,-1){12}}

\put(45,20){\line(2,1){12}}

\put(45,20){\line(2,-1){12}}

\put(18,14){\line(0,1){12}}

\put(57,14){\line(0,1){12}}

\put(27,6){{\scriptsize {\bf {Fig.7}}: $G_{1}$}}

\end{picture}
}
\put(47,-11){\begin{picture}(50,20)

\put(18,26){\circle*{1.5}}

\put(18,14){\circle*{1.5}}

\put(30,20){\circle*{1.5}}

\put(45,20){\circle*{1.5}}

\put(57,26){\circle*{1.5}}

\put(57,14){\circle*{1.5}}

\qbezier(18,26)(35,30)(45,20)

\put(30,20){\line(1,0){15}}

\put(30,20){\line(-2,1){12}}

\put(30,20){\line(-2,-1){12}}

\put(45,20){\line(2,1){12}}

\put(45,20){\line(2,-1){12}}

\put(18,14){\line(0,1){12}}

\put(57,14){\line(0,1){12}}

\put(27,6){{\scriptsize {\bf {Fig.8}}: G}}

\end{picture}
}

\end{picture}
\end{center}

\end{footnotesize}

In this paper, we will do some research on the above two problems.
The following is a Lemma which is obtained by Liu
\cite{liu1}\cite{liu2} and Xuong \cite{xuo} independently.

\medskip

{\bf Lemma 1.1}  \ \ \ Let \emph{G} be a connected graph, then

 1) $\gamma_{M}( G )$ = $\frac{\beta(G) - \xi(G)}{2}$;

 2) \emph{G} is upper embeddable if and only if  $\xi( \emph{G}
 ) \leqslant 1$, or $G$ has a splitting tree.

 \bigskip

\noindent {\bf 2. Vertex splitting and upper embeddability}

\bigskip
As described in the introduction, an upper embeddable graph may be
changed into a non-upper embeddable graph after a vertex splitting
operation. How does a graph remain the upper embeddability after
vertex splitting operations? In this section, we provide some
results on this problem.

\medskip

{\bf Lemma  2.1}   \ \ Let $G$ be an upper embeddable graph, $v$ be
a vertex of $G$ with deg$_G(v)\geqslant$3, and $v_1, v_2, \dots,
v_{n}$ be all the neighbors of $v$ in $G$. If the $v$-$local$
subgraph $G_{loc}(v)$ is connected, then there must exist a
splitting tree $\mathbb{T}$ of $G$ such that all of \{$vv_1, vv_2,
\dots, vv_{n}$\} are edges of $\mathbb{T}$.

\medskip
{\bf Proof } \ \ Let $T$ be an arbitrary splitting tree of $G$.
Since $v_1, v_2, \dots, v_{n}$ are all the neighbors of $v$ in $G$,
the splitting tree $T$ must contain at least one of $\{vv_{i}|i=1,
2, \dots,n\}$ as its edge. Without loss of generality, it may be
assumed that $vv_1\in E(T)$.

If each of $\{vv_{i}|i=2, \dots,n\}$ is an edge of $T$, then the
splitting tree $T$ is $\mathbb{T}$ itself.

If some edges of $\{vv_{i}|i=2, \dots,n\}$ are not in $T$, then
assume, without loss of generality, that $vv_{i_1}, vv_{i_2}, \dots,
vv_{i_m}(m\leqslant n-1)$ are all the edges of $\{vv_{i}|i=2,
\dots,n\}$ which are not in $T$, where the vertex set \{$v_{i_1},
v_{i_2}, \dots, v_{i_m}$\}$\subseteq$ \{$v_2, \dots, v_{n}$\}. Let
$v_{i_j}$ be an arbitrary vertex of  \{$v_{i_1}, v_{i_2}, \dots,
v_{i_m}$\}. Because there is exactly one $u,\omega$-$path$ in $T$
for any two vertices $u$ and $\omega$ in $G$, and the edge
$vv_{i_j}$ is not in $T$, there must be a $vv_{i_j}$-$path$ in $T$,
and the $vv_{i_j}$-$path$ in $T$ must be the style: $v\dots
v_{\alpha}v_{i_j}$, where $v_{\alpha}$ is a vertex of $\{V(G)-\{v,
v_{i_j}\}$\}. Let $T_{i_j}=\{T-v_{\alpha}v_{i_j}\}\cup vv_{i_j}$. It
is obvious that $T_{i_j}$ is a spanning tree of $G$ and the edge
$vv_{i_j}\in E(T_{i_j})$. Through series of processes similar to
that of getting $T_{i_j}$, a spanning tree $T^{*}$ is obtained,
where all of $\{vv_{1}, vv_2, \dots, vv_{n}\}$ are edges of $T^{*}$.
Since all edges of $\{vv_{1}, vv_2, \dots, vv_{n}\}$ are in $T^{*}$,
each edge of $G_{loc}(v)$ is not in $T^{*}$, or else the spanning
tree $T^{*}$ will contain cycles. So all edges of $G_{loc}(v)$ are
co-tree edges of $T^{*}$. Because the $v$-$local$ subgraph
$G_{loc}(v)$ is connected, we can get that $\xi(G,T^{*})\leqslant
\xi(G,T) = \xi(G) \leqslant 1$. So $T^{*}$ is a splitting tree of
$G$ which satisfies the Lemma.$\hspace*{\fill} \Box$

\medskip

{\bf Lemma  2.2}   \ \ Let $G$ be an upper embeddable graph with
minimum degree at least 3, $v$ be a vertex of $G$ with
deg$_{G}(v)$=4, $G^{*}$ be the graph obtained from $G$ by splitting
$v$ into two adjacent vertices $v'$ and $v''$. If the
$splitting$-$edge$ $v'v''$ is not a cut-edge of the $v$-$splitting$
subgraph $G^{*}_{spl}(v)$, then $G^{*}$ is upper embeddable.

\medskip

{\bf Proof } \ \ Let $v_{1}$, $v_{2}$, $v_{3}$, $v_{4}$ be the four
vertices adjacent to $v$ in $G$, and $\mathbb{T}$ be a splitting
tree of $G$. Since $v'v''$ is not a cut-edge of the $v$-$splitting$
subgraph $G^{*}_{spl}(v)$, $G^{*}_{spl}(v)$ must contain at least
one cycle which has $v'v''$ as one of its edges. Without loss of
generality, let $v_{i_1}v_{i_2}v''v'$ be the 4-cycle of
$G^{*}_{spl}(v)$, which is depicted, for example, in Fig.9 or
Fig.11, where \{$v_{i_1}, v_{i_2}$\}=\{$v_1, v_2$\}. Because $G^{*}$
is obtained from $G$ through vertex splitting operation on $v$,
$v_1v_2v$ must be a 3-cycle of $G$, which is depicted, for example,
in Fig.10. In graph $G$, let $C_i$$(i=1,2,3,4)$ denote the connected
component which is obtained from such connected component of
$G-E(\mathbb{T})$ that contains $v_{i}$ as one of its vertices, by
deleting the edges $vv_1, vv_2, vv_3, vv_4, v_1v_2$ from it. It is
possible that $C_{i}$ and $C_{j}$ may be the same connected
component of $G-E(\mathbb{T})$ $(i, j=1,2,3,4 \ $and$ \ i\neq j)$.
If $G$ is upper embeddable, the graph $G^{*}$ in Fig.11, which is
obtained from $G$ through vertex splitting on $v$, is upper
embeddable, for $G^{*}$ can also be viewed as a subdivision of $G$.
So, we should only discuss the upper embeddability of $G^{*}$ in
Fig.9. For $v_{1}$, $v_{2}$, $v_{3}$, $v_{4}$ being all the
neighbors of $v$ in graph $G$, the splitting tree $\mathbb{T}$ of
$G$ must contain at least one edge which belongs to the edge set
$E(v)$=\{$vv_{i}|i=1,2,3,4$\}. It will be discussed in three cases
according to whether at least three edges of $E(v)$ are in
$\mathbb{T}$, or exactly two edges of $E(v)$ are in $\mathbb{T}$, or
only one edge of $E(v)$ is in $\mathbb{T}$. Without loss of
generality, let the edges $v'v_{i_1}$, $v''v_{i_2}$, $v''v_{3}$,
$v'v_{4}$ in $G^{*}$ be the replacement of $vv_{1}$, $vv_{2}$,
$vv_{3}$, $vv_{4}$ in $G$ after vertex splitting on $v$, where the
edge set $\{v'v_{i_1}, v''v_{i_2}\}$ may be $\{v'v_{1}, v''v_{2}\}$
or $\{v'v_{2}, v''v_{1}\}$.

\begin{footnotesize}
\setlength{\unitlength}{0.6mm}
\begin{center}
\begin{picture}(70,34)

\put(-75,-3){\begin{picture}(30,40)

\put(16,28){\circle*{1.5}}

  \put(17,24){$v_{i_1}$}

   \put(6.3,30.3){$C_{i_1}$}

\multiput(16,28)(-0.7,1.8){6}{\circle*{0.5}} 

\put(16,12){\circle*{1.5}}

  \put(17,14.5){$v_{i_2}$}

  \put(10.5,5.5){$C_{i_2}$}

\multiput(16,12)(-1.7,-0.4){6}{\circle*{0.5}}  

\put(30,28){\circle*{1.5}}

  \put(28.5,30){$v'$}

\put(30,12){\circle*{1.5}}

  \put(29,6.5){$v''$}

\multiput(30,28)(-1.8,0){9}{\circle*{0.5}}  

\multiput(30,12)(-1.8,0){9}{\circle*{0.5}}   

\multiput(30,12)(0,1.9){9}{\circle*{0.5}}   

\put(30,28){\line(1,0){14}}  

\put(30,12){\line(1,0){14}}   

\put(44,28){\circle*{1.5}}

   \put(38.5,24){$v_{4}$}

\put(46,37){\circle*{1.5}}

   \put(47.5,37){$v_{l}$}

    \put(48,28){$C_4$}

\multiput(44,28)(0.4,1.9){5}{\circle*{0.5}}  

\multiput(46,37)(-1.7,0.5){5}{\circle*{0.5}}  

\put(44,12){\circle*{1.5}}

  \put(38.5,14){$v_{3}$}

\put(54,11){\circle*{1.5}}

  \put(56,9){$v_{p}$}

   \put(45.5,5){$C_{3}$}

\multiput(44,12)(1.6,-0.1){6}{\circle*{0.5}} 

\multiput(54,11)(0.5,1.5){6}{\circle*{0.5}}    

\put(8.5,35.5){\oval(21,20)[br]}

\put(8.5,4.5){\oval(21,20)[tr]}

\put(51.5,35.5){\oval(21,20)[bl]}

\put(51.5,4.5){\oval(21,20)[tl]}

\put(16,12){\line(0,1){16}} 

\put(20,-4){\scriptsize \bf {Fig.9}: $G^{*}$}

\end{picture}
}

\put(4,-3){\begin{picture}(27,40)

\put(16,28){\circle*{1.5}}

  \put(19,27.5){$v_{1}$}

   \put(8.5,31){$C_{1}$}

\multiput(16,28)(-1.8,0.2){6}{\circle*{0.5}} 

\put(16,12){\circle*{1.5}}

  \put(19,10){$v_{2}$}

    \put(7,8){$C_{2}$}

\multiput(16,12)(-0.7,-1.9){5}{\circle*{0.5}} 

\put(30,20){\circle*{1.5}}

  \put(29,22){$v$}

\put(44,28){\circle*{1.5}}

   \put(37,28.5){$v_{4}$}

\put(46,37){\circle*{1.5}}

   \put(47.5,37){$v_{l}$}

    \put(48,28){$C_4$}

\multiput(44,28)(0.5,1.9){5}{\circle*{0.5}}  

\multiput(46,37)(-1.8,0.5){5}{\circle*{0.5}}  

\put(44,12){\circle*{1.5}}

  \put(37,10.5){$v_{3}$}

\put(54,11){\circle*{1.5}}

  \put(56,9){$v_{p}$}

   \put(46,4){$C_{3}$}

\multiput(44,12)(1.7,-0.2){6}{\circle*{0.5}} 

\multiput(54,11)(0.7,1.7){5}{\circle*{0.5}}    

\put(8.5,35.5){\oval(21,20)[br]}

\put(8.5,4.5){\oval(21,20)[tr]}

\put(51.5,35.5){\oval(21,20)[bl]}

\put(51.5,4.5){\oval(21,20)[tl]}

\multiput(30,20)(-1.7,1){8}{\circle*{0.5}}  

\multiput(30,20)(-1.7,-1){8}{\circle*{0.5}}   

\put(30,20){\line(5,3){14}}    

\put(30,20){\line(5,-3){14}}   

\put(16,12){\line(0,1){16}}    

\put(63,20){\vector(1,0){12}}

\put(0,20){\vector(-1,0){12}}

\put(22,-4){\scriptsize \bf {Fig.10}: $G$}

\end{picture}
}
\put(80,-3){\begin{picture}(27,40)

\put(16,28){\circle*{1.5}}

  \put(17,24){$v_{i_1}$}

\multiput(16,28)(-1.7,1.5){5}{\circle*{0.5}}    

\put(16,12){\circle*{1.5}}

  \put(17,14.5){$v_{i_2}$}

\multiput(16,12)(-1.7,-1.5){5}{\circle*{0.5}} 

\put(30,28){\circle*{1.5}}

  \put(28.5,30){$v'$}

\put(30,12){\circle*{1.5}}

  \put(29,6.5){$v''$}

\multiput(30,28)(-1.8,0){9}{\circle*{0.5}}    

\multiput(30,12)(-1.8,0){9}{\circle*{0.5}}    

\multiput(30,12)(0,1.9){9}{\circle*{0.5}}    

\multiput(30,28)(1.8,0){9}{\circle*{0.5}}    

\qbezier(30,28)(43,20)(44,12)                 

\put(44,28){\circle*{1.5}}

   \put(38.5,24){$v_{4}$}

\put(51,35){\circle*{1.5}}

   \put(53,32){$v_{l}$}

\put(44,28){\line(1,1){7}}   

\multiput(51,35)(-1.8,1.3){5}{\circle*{0.5}}    

\put(44,12){\circle*{1.5}}

  \put(37,13){$v_{3}$}

\put(53,8){\circle*{1.5}}

  \put(52,4){$v_{p}$}

\multiput(44,12)(1.7,-0.9){5} {\circle*{0.5}}    

\multiput(53,8)(1,1.5){6}{\circle*{0.5}}    

\put(8.5,35.5){\oval(21,20)[br]}

\put(8.5,4.5){\oval(21,20)[tr]}

\put(51.5,35.5){\oval(21,20)[bl]}

\put(51.5,4.5){\oval(21,20)[tl]}

\put(16,12){\line(0,1){16}}

\put(22,-4){\scriptsize \bf {Fig.11}: $G^{*}$}

\end{picture}
}

\end{picture}
\end{center}

\end{footnotesize}

\medskip

{\bf Case 1:} \ At least three edges of $E(v)$ are in $\mathbb{T}$.

\medskip
Without loss of generality, let $vv_{1}$, $vv_{2}$, $\dots$,
$vv_{n}$($n=$3 or 4) be all the edges of $E(v)$ which are in
$\mathbb{T}$. Obviously, if exactly three edges of $E(v)$, which are
denoted by $E_3(v)$, are in $\mathbb{T}$, and $E^{*}_3(v)$ denotes
the replacement of $E_3(v)$ after vertex splitting on $v$ in $G$,
then $T^{*} = (G^{*} \cap \mathbb{T}) \cup v'v'' \cup E^{*}_3(v)$ is
a spanning tree of $G^{*}$. If the four edges of $E(v)$ are all in
$\mathbb{T}$, $T^{*} = (G^{*} \cap \mathbb{T}) \cup v'v'' \cup
\{v'v_{i_1}, v''v_{i_2}, v''v_{3}, v'v_4\}$ is a spanning tree of
$G^{*}$. Furthermore, $\xi(G^{*},T^{*})$ = $\xi(G,\mathbb{T}) =
\xi(G) \leqslant 1$. So $T^{*}$ is a splitting tree of $G^{*}$, and
in Case 1 $G^{*}$ is upper embeddable.

\medskip

{\bf Case 2:} \  Exactly two edges of $E(v)$ are in $\mathbb{T}$.

\medskip

The two edges of $E(v)$ in $\mathbb{T}$ may be (i) $vv_1$ and
$vv_2$; or (ii) $vv_3$ and $vv_4$; or (iii) one edge belongs to
\{$vv_1, vv_2$\} and the other belongs to \{$vv_3, vv_4$\}.

\medskip

\textsf{Subcase 2.1:} \ The two edges of $E(v)$ in $\mathbb{T}$ are
$vv_1$ and $vv_2$.

\medskip

In this case, the edge $v_1v_2$ in $G$ can not be an edge of
$\mathbb{T}$, or else $vv_1v_2$ would form a 3-cycle of
$\mathbb{T}$. Let $G^{*}$, which is depicted in Fig.9, denotes the
graph obtained from $G$ through vertex splitting on $v$, where
\{$C_{i_1}, C_{i_2}$\}=\{$C_{1}, C_{2}$\}, and \{$v_{i_1},
v_{i_2}$\}=\{$v_{1}, v_{2}$\}.
\medskip

\textsf{\textit{Subcase 2.1.1:}} \ $C_{3}$ and  $C_{4}$ are the same
connected component of $G$.

\medskip

In this case, let $T^{*} = (G^{*} \cap \mathbb{T}) \cup v'v'' \cup
\{v'v_{i_1}, v''v_{i_2}\}$. It is obvious that $T^{*}$ is a spanning
tree of $G^{*}$, and $\xi(G^{*},T^{*})$ = $\xi(G,\mathbb{T}) =
\xi(G) \leqslant 1$. So $T^{*}$ is a splitting tree of $G^{*}$, and
$G^{*}$ is upper embeddable in Subcase-2.1.1.
\medskip

\textsf{\textit{Subcase 2.1.2:}} \ $C_{3}$ and  $C_{4}$ are two
different connected components of $G$.
\medskip

In graph $G^{*}$, if at least one of  $C_3\cup v''v_3$ and  $C_4\cup
v'v_4$ contains an even number of edges, then let $T^{*} = (G^{*}
\cap \mathbb{T}) \cup v'v'' \cup \{v'v_{i_1}, v''v_{i_2}\}$. It is
obvious that $\xi(G^{*},T^{*})$ = $\xi(G,\mathbb{T}) = \xi(G)
\leqslant 1$. So $T^{*}$ is a splitting tree of $G^{*}$, and $G^{*}$
is upper embeddable.

If both $C_3\cup v''v_3$ and $C_4\cup v'v_4$ contain an odd number
of edges, then $C_3$ and $C_4$ both contain an even number of edges.
Because there is exactly one $u,\omega$-$path$ in $\mathbb{T}$ for
any two vertices $u$ and $\omega$ in $G$, and both $vv_3$ and $vv_4$
are not in $\mathbb{T}$, there must be exactly one $v,v_3$-$path$ in
$\mathbb{T}$, and the $v,v_3$-$path$ in $\mathbb{T}$ must be of the
form as $vv_1\dots v_{p}v_3$ or $vv_2\dots v_{p}v_3$.  Also, there
must be exactly one $v,v_4$-$path$ in $\mathbb{T}$, and the
$v,v_4$-$path$ in $\mathbb{T}$ must be of the form as $vv_1\dots
v_{l}v_4$ or $vv_2\dots v_{l}v_4$. Furthermore, the $v,v_3$-$path$
and $v,v_4$-$path$ in $\mathbb{T}$ can not form a cycle. It is
discussed in the following three subcases.

\texttt{\textit{Subcase 2.1.2-a:}} \ The $v,v_3$-$path$ and
$v,v_4$-$path$ in $\mathbb{T}$ are $vv_1\dots v_{p}v_3$ and
$vv_1\dots v_{l}v_4$ respectively.

If the edges $vv_1$ and $vv_2$ in $G$ are replaced, after the vertex
splitting on $v$, by $v'v_{i_1}$ and $v''v_{i_2}$ respectively, then
$T^{*}_{1} = (G^{*} \cap \mathbb{T}) \cup \{v'v_4, v''v_3,
v''v_{i_2}\}$ is a spanning tree of $G^{*}$. Noticing that the size
of $C_{i_1}\cup v_{i_1}v_{i_2}\cup C_{i_2}\cup v_{i_1}v'\cup v'v''$
and $C_{i_1}\cup v_{i_1}v_{i_2}\cup C_{i_2}$ have the same parity,
and both the size of $C_3$ and $C_4$ are an even number, we can
easily get that $\xi(G^{*},T^{*}_{1})$ = $\xi(G,\mathbb{T}) = \xi(G)
\leqslant 1$. So $T^{*}_{1}$ is a splitting tree of $G^{*}$, and
$G^{*}$ is upper embeddable.

After the vertex splitting on $v$ in $G$, if the edge $vv_1$ is
replaced by $v''v_{i_2}$, and $vv_2$ by $v'v_{i_1}$ respectively,
then $T^{*}_{2} = (G^{*} \cap \mathbb{T}) \cup \{v'v_4, v'v_{i_1},
v''v_{3}\}$ is a spanning tree of $G^{*}$. It is obvious that
$\xi(G^{*},T^{*}_{2})$ = $\xi(G,\mathbb{T}) = \xi(G) \leqslant 1$.
So $T^{*}_{2}$ is a splitting tree of $G^{*}$, and $G^{*}$ is upper
embeddable.

\texttt{\textit{Subcase 2.1.2-b:}} \ The $v,v_3$-$path$ and
$v,v_4$-$path$ in $\mathbb{T}$ are $vv_2\dots v_{p}v_3$ and
$vv_1\dots v_{l}v_4$ respectively.

In this case, let $T^{*} = (G^{*} \cap \mathbb{T}) \cup \{v'v_4,
 v'v'', v''v_3\}$ be a spanning tree of $G^{*}$. It is obvious that
$\xi(G^{*},T^{*})$ = $\xi(G,\mathbb{T}) = \xi(G) \leqslant 1$. So
$T^{*}$ is a splitting tree of $G^{*}$, and $G^{*}$ is upper
embeddable.

\texttt{\textit{Subcase 2.1.2-c:}} \ The $v,v_3$-$path$ and
$v,v_4$-$path$ in $\mathbb{T}$ are $vv_2\dots v_{p}v_3$ and
$vv_2\dots v_{l}v_4$ respectively, or $vv_1\dots v_{p}v_3$ and
$vv_2\dots v_{l}v_4$ respectively.

In this case, it is similar to that of Subcase 2.1.2-a and Subcase
2.1.2-b to get that $G^{*}$ contains a splitting tree.

So, in Subcase-2.1.2, $G^{*}$ is upper embeddable.

\medskip

\textsf{Subcase 2.2:} \ The two edges of $E(v)$ in $\mathbb{T}$ are
$vv_3$ and $vv_4$.

\medskip

In this case, according to $v_1v_2$ being an edge of $\mathbb{T}$ or
not, it will be discussed in  the following two subcases.

\medskip

\textsf{\textit{Subcase 2.2.1:}} \ The edge $v_1v_2$ of $G$ is not
in $\mathbb{T}$.

\medskip

In this case, let $T^{*} = (G^{*} \cap \mathbb{T}) \cup \{v'v_4,
v'v'', v''v_3\}$ be a spanning tree of $G^{*}$. It is obvious that
$\xi(G^{*},T^{*})$ = $\xi(G,\mathbb{T}) = \xi(G) \leqslant 1$. So
$T^{*}$ is a splitting tree of $G^{*}$.

\medskip
\textsf{\textit{Subcase 2.2.2:}} \ The edge $v_1v_2$ of $G$ is an
edge of $\mathbb{T}$.

\medskip

It will be discussed in the following subcases.

\texttt{\textit{Subcase 2.2.2-1:}} \  $C_{i_1}$ and $C_{i_2}$ are
the same connected component of $G$.

In this case, let $T^{*} = (G^{*} \cap \mathbb{T}) \cup \{v'v_4,
v'v'', v''v_3\}$ be a spanning tree of $G^{*}$. It is obvious that
$\xi(G^{*},T^{*})$ = $\xi(G,\mathbb{T}) = \xi(G) \leqslant 1$. So
$T^{*}$ is a splitting tree of $G^{*}$.

\texttt{\textit{Subcase 2.2.2-2:}} \  $C_{i_1}$ and $C_{i_2}$ are
two different connected components of $G$.

If at least one of  $C_{i_1}\cup v'v_{i_1}$ and  $C_{i_2}\cup
v''v_{i_2}$ contains an even number of edges, then let $T^{*} =
(G^{*} \cap \mathbb{T}) \cup \{v'v_4, v'v'', v''v_3\}$. It is
obvious that $\xi(G^{*},T^{*})$ = $\xi(G,\mathbb{T}) = \xi(G)
\leqslant 1$. So $T^{*}$ is a splitting tree of $G^{*}$, and $G^{*}$
is upper embeddable.

If both $C_{i_1}\cup v'v_{i_1}$ and  $C_{i_2}\cup v''v_{i_2}$
contain an odd number of edges, then $C_{i_1}$ and $C_{i_2}$ both
contain an even number of edges. Because there is exactly one
$u,\omega$-$path$ in $\mathbb{T}$ for any two vertices $u$ and
$\omega$ in $G$, and both $vv_1$ and $vv_2$ are not in $\mathbb{T}$,
there must be exactly one $v,v_1$-$path$ in $\mathbb{T}$, and this
$v,v_1$-$path$ in $\mathbb{T}$ may be the form as $vv_4\dots
v_{1}v_2$, or $vv_4\dots v_{2}v_1$, or $vv_3\dots v_{1}v_2$, or
$vv_3\dots v_{2}v_1$. It is discussed in the following two subcases.

\texttt{\textit{Subcase 2.2.2-2a:}} \ The $v,v_1$-$path$ in
$\mathbb{T}$ is $vv_4\dots v_{1}v_2$ or $vv_4\dots v_{2}v_1$.

In this case, let $T^{*} = (G^{*} \cap \mathbb{T}) \cup \{v'v_{4},
v''v_3, v''v_{i_2}\}$. Noticing that both $C_{i_1}\cup v_{i_1}v'\cup
v'v''$ and $C_{i_2}$ contain an even number of edges, we can get
that $\xi(G^{*},T^{*})$ = $\xi(G,\mathbb{T}) = \xi(G) \leqslant 1$.
So $T^{*}$ is a splitting tree of $G^{*}$, and $G^{*}$ is upper
embeddable.

\texttt{\textit{Subcase 2.2.2-2b:}} \ The $v,v_1$-$path$ in
$\mathbb{T}$ is $vv_3\dots v_{1}v_2$ or $vv_3\dots v_{2}v_1$.

In this case, let $T^{*} = (G^{*} \cap \mathbb{T}) \cup \{v'v_{i_1},
v'v_{4}, v''v_{3}\}$. It is obvious that $\xi(G^{*},T^{*})$ =
$\xi(G,\mathbb{T}) = \xi(G) \leqslant 1$. So $T^{*}$ is a splitting
tree of $G^{*}$, and $G^{*}$ is upper embeddable.

\medskip

\textsf{Subcase 2.3:} \ The two edges of $E(v)$ in $\mathbb{T}$ are
such two edges that one is selected from \{$vv_1, vv_2,$\} and the
other is selected from \{$vv_3, vv_4$\}.

\medskip

Without loss of generality, let the two edges of $E(v)$ in
$\mathbb{T}$ are $vv_1$ and $vv_3$, which is illustrated in Fig.13.
We will discuss in the following two subcases.

\medskip

\textsf{\textit{Subcase 2.3.1:}} \ After the vertex splitting on $v$
in $G$, the replacements of $vv_1$ and $vv_3$ are both adjacent to
$v'$ or both adjacent to $v''$.

\medskip

Without loss of generality, let the replacements of $vv_1$ and
$vv_3$ are both adjacent to $v''$, which is illustrated in Fig.12.
Let $T^{*} = (G^{*} \cap \mathbb{T}) \cup \{v''v_3, v''v_{i_2},
v'v''\}$. It is obvious that $\xi(G^{*},T^{*})$ = $\xi(G,\mathbb{T})
= \xi(G) \leqslant 1$. So $T^{*}$ is a splitting tree of $G^{*}$.

\medskip
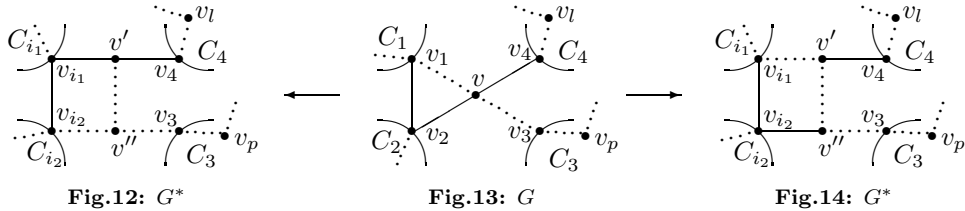
\begin{footnotesize}

\setlength{\unitlength}{0.6mm}
\begin{center}
\begin{picture}(70,31)

\put(-75,-3){\begin{picture}(30,40)

\put(16,28){\circle*{1.5}}

  \put(17,24){$v_{i_1}$}

   \put(6.3,30.3){$C_{i_1}$}

\multiput(16,28)(-0.9,1.8){5}{\circle*{0.5}} 

\put(16,12){\circle*{1.5}}

  \put(17,14.5){$v_{i_2}$}

  \put(10.5,5.5){$C_{i_2}$}

\multiput(16,12)(-2,-0.4){5}{\circle*{0.5}}  

\put(30,28){\circle*{1.5}}

  \put(28.5,30){$v'$}

\put(30,12){\circle*{1.5}}

  \put(29,6.5){$v''$}

\put(30,28){\line(-1,0){14}}  

\multiput(30,12)(-2,0){8}{\circle*{0.5}}   

\multiput(30,12)(0,2){8}{\circle*{0.5}}   

\put(30,28){\line(1,0){14}}  

\multiput(30,12)(2,0){8}{\circle*{0.5}}   

\put(44,28){\circle*{1.5}}

   \put(38.5,24){$v_{4}$}

\put(46,37){\circle*{1.5}}

   \put(47.5,37){$v_{l}$}

    \put(48,28){$C_4$}

\multiput(44,28)(0.6,2.2){4}{\circle*{0.5}}  

\multiput(46,37)(-1.8,0.6){5}{\circle*{0.5}}  

\put(44,12){\circle*{1.5}}

  \put(38.5,14){$v_{3}$}

\put(54,11){\circle*{1.5}}

  \put(56,9){$v_{p}$}

   \put(45.5,5){$C_{3}$}

\multiput(44,12)(1.8,-0.1){6}{\circle*{0.5}} 

\multiput(54,11)(0.6,1.8){5}{\circle*{0.5}}    

\put(8.5,35.5){\oval(21,20)[br]}

\put(8.5,4.5){\oval(21,20)[tr]}

\put(51.5,35.5){\oval(21,20)[bl]}

\put(51.5,4.5){\oval(21,20)[tl]}

\put(16,12){\line(0,1){16}} 

\put(19.5,-4){\scriptsize \bf {Fig.12}: $G^{*}$}

\end{picture}
}

\put(4,-3){\begin{picture}(27,40)

\put(16,28){\circle*{1.5}}

  \put(19,27.5){$v_{1}$}

   \put(8.5,31){$C_{1}$}

\multiput(16,28)(-2,0.2){5}{\circle*{0.5}} 

\put(16,12){\circle*{1.5}}

  \put(19,10){$v_{2}$}

    \put(7,8){$C_{2}$}

\multiput(16,12)(-0.8,-1.8){5}{\circle*{0.5}} 

\put(30,20){\circle*{1.5}}

  \put(29,22){$v$}

\put(44,28){\circle*{1.5}}

   \put(37,28.5){$v_{4}$}

\put(46,37){\circle*{1.5}}

   \put(47.5,37){$v_{l}$}

    \put(48,28){$C_4$}

\multiput(44,28)(0.6,2.2){4}{\circle*{0.5}}  

\multiput(46,37)(-1.8,0.6){5}{\circle*{0.5}}  

\put(44,12){\circle*{1.5}}

  \put(37,10.5){$v_{3}$}

\put(54,11){\circle*{1.5}}

  \put(56,9){$v_{p}$}

   \put(46,4){$C_{3}$}

\multiput(44,12)(1.8,-0.1){6}{\circle*{0.5}} 

\multiput(54,11)(0.6,1.8){5}{\circle*{0.5}}    

\put(8.5,35.5){\oval(21,20)[br]}

\put(8.5,4.5){\oval(21,20)[tr]}

\put(51.5,35.5){\oval(21,20)[bl]}

\put(51.5,4.5){\oval(21,20)[tl]}

\multiput(30,20)(-1.7,1){8}{\circle*{0.5}}  

\put(30,20){\line(-5,-3){14}}   

\put(30,20){\line(5,3){14}}    

\multiput(30,20)(1.7,-1){8}{\circle*{0.5}}   

\put(16,12){\line(0,1){16}}    

\put(63,20){\vector(1,0){12}}

\put(0,20){\vector(-1,0){12}}

\put(19.5,-4){\scriptsize \bf {Fig.13}: $G$}

\end{picture}
}
\put(80,-3){\begin{picture}(30,40)

\put(16,28){\circle*{1.5}}

  \put(17,24){$v_{i_1}$}

   \put(6.3,30.3){$C_{i_1}$}

\multiput(16,28)(-0.8,1.8){5}{\circle*{0.5}} 

\put(16,12){\circle*{1.5}}

  \put(17,14){$v_{i_2}$}

  \put(10.5,5.5){$C_{i_2}$}

\multiput(16,12)(-2,-0.5){5}{\circle*{0.5}}  

\put(30,28){\circle*{1.5}}

  \put(28.5,30){$v'$}

\put(30,12){\circle*{1.5}}

  \put(29,6.5){$v''$}

\multiput(30,28)(-2,0){8}{\circle*{0.5}}  

\put(30,12){\line(-1,0){14}}   

\multiput(30,12)(0,2){9}{\circle*{0.5}}   

\put(30,28){\line(1,0){14}}  

\multiput(30,12)(2,0){8}{\circle*{0.5}}   

\put(44,28){\circle*{1.5}}

   \put(38.5,24){$v_{4}$}

\put(46,37){\circle*{1.5}}

   \put(47.5,37){$v_{l}$}

    \put(48,28){$C_4$}

\multiput(44,28)(0.6,2.2){4}{\circle*{0.5}}  

\multiput(46,37)(-1.8,0.6){5}{\circle*{0.5}}  

\put(44,12){\circle*{1.5}}

  \put(38.5,14){$v_{3}$}

\put(54,11){\circle*{1.5}}

  \put(56,9){$v_{p}$}

   \put(45.5,5){$C_{3}$}

\multiput(44,12)(1.8,-0.1){6}{\circle*{0.5}} 

\multiput(54,11)(0.6,1.8){5}{\circle*{0.5}}    

\put(8.5,35.5){\oval(21,20)[br]}

\put(8.5,4.5){\oval(21,20)[tr]}

\put(51.5,35.5){\oval(21,20)[bl]}

\put(51.5,4.5){\oval(21,20)[tl]}

\put(16,12){\line(0,1){16}} 

\put(19.5,-4){\scriptsize \bf {Fig.14}: $G^{*}$}

\end{picture}
}

\end{picture}
\end{center}

\end{footnotesize}

\medskip

\textsf{\textit{Subcase 2.3.2:}} \ After the vertex splitting on $v$
in $G$, the replacements of $vv_1$ and $vv_3$ are  adjacent to $v'$
and $v''$ respectively.
\medskip

Without loss of generality, let $vv_1$ and $vv_3$ be replaced, after
vertex splitting on $v$, by $v'v_{i_1}$ and $v''v_3$ respectively,
which is illustrated in Fig.14.

\texttt{\textit{Subcase 2.3.2-1:}} \ In graph $G$, the edge $v_1v_2$
is not an edge of $\mathbb{T}$.

If $C_{4}$ and one of \{$C_{i_1}, C_{i_2}$\} are the same connected
component of $G$, then $T^{*}_{1} = (G^{*} \cap \mathbb{T}) \cup
\{v'v_{i_1}, v'v'', v''v_3\}$ is a splitting tree of $G^{*}$.

If $C_{4}$ is a connected component of $G$ which is different from
both of \{$C_{i_1}, C_{i_2}$\}, we will discuss in two subcases.

\texttt{\textit{Subcase 2.3.2-1a:}} \ At least one of $C_{i_1}\cup
v_{i_1}v_{i_2}\cup C_{i_2}\cup v_{i_2}v''$ and  $C_{4}\cup v'v_{4}$
contains an even number of edges.

In this case, let $T^{*} = (G^{*} \cap \mathbb{T}) \cup \{v'v_{i_1},
v'v'', v''v_3\}$. It is obvious that $\xi(G^{*},T^{*})$ =
$\xi(G,\mathbb{T}) = \xi(G) \leqslant 1$. So $T^{*}$ is a splitting
tree of $G^{*}$, and $G^{*}$ is upper embeddable.

\texttt{\textit{Subcase 2.3.2-1b:}} \  Both $C_{i_1}\cup
v_{i_1}v_{i_2}\cup C_{i_2}\cup v_{i_2}v''$ and  $C_{4}\cup v'v_{4}$
contain an odd number of edges.

In this case, $C_{4}$ contains an even number of edges. Because
there is exactly one $u,\omega$-$path$ in $\mathbb{T}$ for any two
vertices $u$ and $\omega$ in $G$, and both $vv_2$ and $vv_4$ are not
in $\mathbb{T}$, there must be exactly one $v,v_4$-$path$ in
$\mathbb{T}$, and the $v,v_4$-$path$ in $\mathbb{T}$ must be the
form as $vv_1\dots v_{4}$ or $vv_3\dots v_4$. If the $v,v_4$-$path$
in $\mathbb{T}$ is $vv_1\dots v_{4}$, then $T^{*}_{1} = (G^{*} \cap
\mathbb{T}) \cup \{v'v_{4}, v'v'', v''v_3\}$ is a splitting tree of
$G^{*}$. If the $v,v_4$-$path$ in $\mathbb{T}$ is $vv_3\dots v_{4}$,
then $T^{*}_{2} = (G^{*} \cap \mathbb{T}) \cup \{v''v_3, v'v_{4},
v'v_{i_1}\}$ is a splitting tree of $G^{*}$.

\texttt{\textit{Subcase 2.3.2-2:}} \ In graph $G$, the edge $v_1v_2$
is an edge of $\mathbb{T}$.

If at least one of  $C_{i_2}\cup v''v_{i_2}$ and  $C_{4}\cup
v'v_{4}$ contains an even number of edges, then let $T^{*}_{1} =
(G^{*} \cap \mathbb{T}) \cup \{v'v_{i_1}, v'v'', v''v_3\}$. It is
obvious that $\xi(G^{*},T^{*}_{1})$ = $\xi(G,\mathbb{T}) = \xi(G)
\leqslant 1$. So $T^{*}_{1}$ is a splitting tree of $G^{*}$, and
$G^{*}$ is upper embeddable.

If both $C_{i_2}\cup v''v_{i_2}$ and  $C_{4}\cup v'v_{4}$ contain an
odd number of edges, then $C_{i_2}$ and $C_{4}$ both contain an even
number of edges. Let $T^{*}_{2} = (G^{*} \cap \mathbb{T}) \cup
\{v'v_{i_1}, v''v_{i_2}, v''v_3\}$. It is obvious that
$\xi(G^{*},T^{*}_{2})$ = $\xi(G,\mathbb{T}) = \xi(G) \leqslant 1$.
So $T^{*}_{2}$ is a splitting tree of $G^{*}$, and $G^{*}$ is upper
embeddable.

\medskip

{\bf Case 3:}  \ Only one edge of $E(v)$ is in $\mathbb{T}$.

\medskip

According to this edge is selected from \{$vv_1, vv_2$\} or \{$vv_3,
vv_4$\}, it will be discussed in the following Subcase-3.1 and
Subcase-3.2.

\medskip

\textsf{Subcase 3.1:} \ One of \{$vv_1, vv_2$\} is the edge in
$\mathbb{T}$.

\medskip

Without loss of generality, let $vv_1$ be the edge in $\mathbb{T}$,
which is depicted in Fig.16. In addition, throughout Subcase 3.1,
let $vv_1$ and $vv_2$ be replaced by $v'v_{1}$ and $v''v_{2}$
respectively after the vertex splitting on $v$ in $G$; and the edge
set \{$vv_3, vv_4$\} be replaced by \{$v''v_{i_3}, v'v_{i_4}$\},
where \{$v_{i_3}, v_{i_4}$\}=\{$v_{3}, v_{4}$\} and \{$C_{i_3},
C_{i_4}$\}=\{$C_{3}, C_{4}$\}, which is depicted in Fig.15.
According to the edge $v_1v_2$ of $G$ is in the splitting tree
$\mathbb{T}$ or not, it will be discussed in the following two
subcases.

\begin{footnotesize}

\medskip
\setlength{\unitlength}{0.6mm}
\begin{center}
\begin{picture}(70,31)

\put(-100,-3){\begin{picture}(30,40)

\put(16,28){\circle*{1.5}}

  \put(17,24){$v_{1}$}

   \put(6.3,30.3){$C_{1}$}

\multiput(16,28)(-0.8,1.8){5}{\circle*{0.5}} 

\put(16,12){\circle*{1.5}}

  \put(17,14){$v_{2}$}

  \put(10.5,5.5){$C_{2}$}

\multiput(16,12)(-2,-0.5){5}{\circle*{0.5}}  

\put(30,28){\circle*{1.5}}

  \put(28.5,30){$v'$}

\put(30,12){\circle*{1.5}}

  \put(29,6.5){$v''$}

\multiput(30,28)(-2,0){8}{\circle*{0.5}}  

\put(30,12){\line(-1,0){14}}   

\multiput(30,12)(0,2){9}{\circle*{0.5}}   

\put(30,28){\line(1,0){14}}  

\put(30,12){\line(1,0){14}}   

\put(44,28){\circle*{1.5}}

   \put(38.5,24){$v_{i_4}$}

\put(46,37){\circle*{1.5}}

   \put(47.5,37){$v_{l}$}

    \put(48,28){$C_{i_4}$}

\multiput(44,28)(0.6,2.2){4}{\circle*{0.5}}  

\multiput(46,37)(-1.8,0.6){5}{\circle*{0.5}}  

\put(44,12){\circle*{1.5}}

  \put(38,14.9){$v_{i_3}$}

\put(54,11){\circle*{1.5}}

  \put(56,9){$v_{p}$}

   \put(45.5,5){$C_{i_3}$}

\multiput(44,12)(1.8,-0.1){6}{\circle*{0.5}} 

\multiput(54,11)(0.6,1.8){5}{\circle*{0.5}}    

\put(8.5,35.5){\oval(21,20)[br]}

\put(8.5,4.5){\oval(21,20)[tr]}

\put(51.5,35.5){\oval(21,20)[bl]}

\put(51.5,4.5){\oval(21,20)[tl]}

\put(16,12){\line(0,1){16}} 

\put(19.5,-4){\scriptsize \bf {Fig.15}: $G^{*}$}

\end{picture}
}

\put(-30,-3){\begin{picture}(27,40)

\put(16,28){\circle*{1.5}}

  \put(19,27.5){$v_{1}$}

   \put(8.5,31){$C_{1}$}

\multiput(16,28)(-2,0.2){6}{\circle*{0.5}} 

\put(16,12){\circle*{1.5}}

  \put(19,10){$v_{2}$}

    \put(7,8){$C_{2}$}

\multiput(16,12)(-0.8,-1.8){6}{\circle*{0.5}} 

\put(30,20){\circle*{1.5}}

  \put(29,22){$v$}

\put(44,28){\circle*{1.5}}

   \put(37,28.5){$v_{4}$}

\put(46,37){\circle*{1.5}}

   \put(47.5,37){$v_{l}$}

    \put(48,28){$C_4$}

\multiput(44,28)(0.6,2.2){4}{\circle*{0.5}}  

\multiput(46,37)(-1.8,0.6){5}{\circle*{0.5}}  

\put(44,12){\circle*{1.5}}

  \put(37,10.5){$v_{3}$}

\put(54,11){\circle*{1.5}}

  \put(56,9){$v_{p}$}

   \put(46,4){$C_{3}$}

\multiput(44,12)(1.8,-0.1){6}{\circle*{0.5}} 

\multiput(54,11)(0.6,1.8){5}{\circle*{0.5}}    

\put(8.5,35.5){\oval(21,20)[br]}

\put(8.5,4.5){\oval(21,20)[tr]}

\put(51.5,35.5){\oval(21,20)[bl]}

\put(51.5,4.5){\oval(21,20)[tl]}

\multiput(30,20)(-1.7,1){8}{\circle*{0.5}}  

\put(30,20){\line(-5,-3){14}}   

\put(30,20){\line(5,3){14}}    

\put(30,20){\line(5,-3){14}}  

\put(16,12){\line(0,1){16}}    

\put(1,20){\vector(-1,0){9}}

\put(19.5,-4){\scriptsize \bf {Fig.16}: $G$}

\end{picture}
}
\put(35,-3){\begin{picture}(30,40)

\put(16,28){\circle*{1.5}}

  \put(17,24){$v_{1}$}

   \put(6.3,30.3){$C_{1}$}

\multiput(16,28)(-0.8,1.8){5}{\circle*{0.5}} 

\put(16,12){\circle*{1.5}}

  \put(17,14){$v_{2}$}

  \put(10.5,5.5){$C_{2}$}

\multiput(16,12)(-2,-0.5){5}{\circle*{0.5}}  

\put(30,28){\circle*{1.5}}

  \put(28.5,30){$v'$}

\put(30,12){\circle*{1.5}}

  \put(29,6.5){$v''$}

\multiput(30,28)(-2,0){8}{\circle*{0.5}}  

\multiput(30,12)(-2,0){8}{\circle*{0.5}}   

\put(30,12){\line(0,1){16}}    

\put(30,28){\line(1,0){14}}  

\put(30,12){\line(1,0){14}}   

\put(44,28){\circle*{1.5}}

   \put(38.5,24){$v_{i_4}$}

\put(46,37){\circle*{1.5}}

   \put(47.5,37){$v_{l}$}

    \put(48,28){$C_{i_4}$}

\multiput(44,28)(0.6,2.2){4}{\circle*{0.5}}  

\multiput(46,37)(-1.8,0.6){5}{\circle*{0.5}}  

\put(44,12){\circle*{1.5}}

  \put(38.2,14.9){$v_{i_3}$}

\put(54,11){\circle*{1.5}}

  \put(56,9){$v_{p}$}

   \put(45.5,5){$C_{i_3}$}

\multiput(44,12)(1.8,-0.1){6}{\circle*{0.5}} 

\multiput(54,11)(0.6,1.8){5}{\circle*{0.5}}    

\put(8.5,35.5){\oval(21,20)[br]}

\put(8.5,4.5){\oval(21,20)[tr]}

\put(51.5,35.5){\oval(21,20)[bl]}

\put(51.5,4.5){\oval(21,20)[tl]}

\multiput(16,12)(0,2){9}{\circle*{0.5}} 

\put(19.5,-4){\scriptsize \bf {Fig.17}: $G^{*}$}

\end{picture}
}

\put(100.3,-3){\begin{picture}(30,40)

\put(16,28){\circle*{1.5}}

  \put(17,24){$v_{1}$}

   \put(6.3,30.3){$C_{1}$}

\multiput(16,28)(-0.8,1.8){5}{\circle*{0.5}} 

\put(16,12){\circle*{1.5}}

  \put(17,14){$v_{2}$}

  \put(10.5,5.5){$C_{2}$}

\multiput(16,12)(-2,-0.5){5}{\circle*{0.5}}  

\put(30,28){\circle*{1.5}}

  \put(28.5,30){$v'$}

\put(30,12){\circle*{1.5}}

  \put(29,6.5){$v''$}

\multiput(30,28)(-2,0){8}{\circle*{0.5}}  

\put(30,12){\line(-1,0){14}}  

\put(30,12){\line(0,1){16}}    

\put(30,28){\line(1,0){14}}  

\multiput(30,12)(2,0){8}{\circle*{0.5}}    

\put(44,28){\circle*{1.5}}

   \put(38.5,24){$v_{i_4}$}

\put(46,37){\circle*{1.5}}

   \put(47.5,37){$v_{l}$}

    \put(48,28){$C_{i_4}$}

\multiput(44,28)(0.6,2.2){4}{\circle*{0.5}}  

\multiput(46,37)(-1.8,0.6){5}{\circle*{0.5}}  

\put(44,12){\circle*{1.5}}

  \put(38.2,14.9){$v_{i_3}$}

\put(54,11){\circle*{1.5}}

  \put(56,9){$v_{p}$}

   \put(45.5,5){$C_{i_3}$}

\multiput(44,12)(1.8,-0.1){6}{\circle*{0.5}} 

\multiput(54,11)(0.6,1.8){5}{\circle*{0.5}}    

\put(8.5,35.5){\oval(21,20)[br]}

\put(8.5,4.5){\oval(21,20)[tr]}

\put(51.5,35.5){\oval(21,20)[bl]}

\put(51.5,4.5){\oval(21,20)[tl]}

\multiput(16,12)(0,2){9}{\circle*{0.5}} 

\put(19.5,-4){\scriptsize \bf {Fig.18}: $G^{*}$}

\end{picture}
}
\end{picture}
\end{center}
\end{footnotesize}

\medskip

\textsf{\textit{Subcase 3.1.1:}} \ In graph $G$, $v_1v_2$ is not an
edge of $\mathbb{T}$. It is discussed in the following subcases.

\medskip

\texttt{\textit{Subcase 3.1.1-1:}} \ In graph $G^{*}$, $C_{1}\cup
v_{1}v_{2}\cup C_{2}\cup v_{2}v''\cup v''v_{i_3}\cup C_{i_3}\cup
C_{i_4}\cup v'v_{i_4}$ contains an odd number of edges.

In this case, $T^{*} = (G^{*} \cap \mathbb{T}) \cup \{v'v_{1},
v'v''\}$ is a splitting tree of $G^{*}$. So, in Subcase 3.1.1-1,
$G^{*}$ is upper embeddable.

\texttt{\textit{Subcase 3.1.1-2:}} \ In graph $G^{*}$, $C_{1}\cup
v_{1}v_{2}\cup C_{2}\cup v_{2}v''\cup v''v_{i_3}\cup C_{i_3}\cup
C_{i_4}\cup v'v_{i_4}$ contains an even number of edges.

In this case, if $C_{i_4}\cup v'v_{i_4}$ contains an even number of
edges, then $T^{*} = (G^{*} \cap \mathbb{T}) \cup \{v'v_{1},
v'v''\}$ is a splitting tree of $G^{*}$.

If $C_{i_4}\cup v'v_{i_4}$ contains an odd number of edges, then
$C_{1}\cup v_{1}v_{2}\cup C_{2}\cup v_{2}v''\cup v''v_{i_3}\cup
C_{i_3}$ contains an odd number of edges too. It is discussed in the
following two subcases.

\texttt{\textit{Subcase 3.1.1-2a:}} \ In graph $G^{*}$, the
connected component $C_{i_4}$ is the same with at least one of
\{$C_{1}$, $C_{i_3}$\}.

In this case, $T^{*} = (G^{*} \cap \mathbb{T}) \cup \{v'v_{1},
v'v''\}$ is a splitting tree of $G^{*}$.

\texttt{\textit{Subcase 3.1.1-2b:}} \ In graph $G^{*}$, neither of
\{$C_{1}$, $C_{i_3}$\} is the same connected component with
$C_{i_4}$.

Because there is exactly one $u,\omega$-$path$ in $\mathbb{T}$ for
any two vertices $u$ and $\omega$ in $G$, and none of \{$vv_2, vv_3,
vv_4$\} is an edge of $\mathbb{T}$, there must be exactly one
$v,v_3$-$path$ and exactly one $v,v_4$-$path$ in $\mathbb{T}$, and
the $v,v_3$-$path$, $v,v_4$-$path$ in $\mathbb{T}$ must be of the
form as $vv_1\dots v_{3}$ and  $vv_1\dots v_4$ respectively.
Noticing that both $C_{i_4}$ and $v'v_{1}\cup C_{1}\cup
v_{1}v_{2}\cup C_{2}\cup v_{2}v''\cup v''v_{i_3}\cup C_{i_3}$ are
connected component of $G^{*}$ with an even number of edges, we can
easily  get that $T^{*} = (G^{*} \cap \mathbb{T}) \cup \{v'v_{i_4},
v'v''\}$ is a splitting tree of $G^{*}$.

\medskip

\textsf{\textit{Subcase 3.1.2:}} \ In graph $G$, $v_1v_2$ is an edge
of $\mathbb{T}$. It is discussed in the following subcases.
\medskip

In graph $G^{*}$, if $C_{i_4}$ is the same connected component with
at least one of \{$C_{1}, C_{2}, C_{i_3}$\}, then $T^{*} = (G^{*}
\cap \mathbb{T}) \cup \{v'v_{1}, v'v''\}$ is a splitting tree of
$G^{*}$. If any pair of components, which is selected from \{$C_{1},
C_{2}, C_{i_3}, C_{i_4}$\}, is not the same connected component of
$G^{*}$, then it will be discussed in the following two subcases.

\texttt{\textit{Subcase 3.1.2-1:}} \ In graph $G^{*}$, $C_{2}\cup
v_{2}v''\cup v''v_{i_3}\cup C_{i_3}\cup C_{i_4}\cup v'v_{i_4}$
contains an odd number of edges.

Noticing that one of \{$C_{i_4}\cup v'v_{i_4}, \ C_{2}\cup
v_{2}v''\cup v''v_{i_3}\cup C_{i_3}$\} is a connected component of
$G^{*}$ which contains an even number of edges, and the other is one
which contains an odd number of edges, we can easily get that $T^{*}
= (G^{*} \cap \mathbb{T}) \cup \{v'v_{1}, v'v''\}$ is a splitting
tree of $G^{*}$.

\texttt{\textit{Subcase 3.1.2-2:}} \ In graph $G^{*}$, $C_{2}\cup
v_{2}v''\cup v''v_{i_3}\cup C_{i_3}\cup C_{i_4}\cup v'v_{i_4}$
contains an even number of edges.

If both $C_{i_4}\cup v'v_{i_4}$ and $C_{2}\cup v_{2}v''\cup
v''v_{i_3}\cup C_{i_3}$ are connected component of $G^{*}$ which
contain an even number of edges, then it is easy to get that $T^{*}
= (G^{*} \cap \mathbb{T}) \cup \{v'v_{1}, v'v''\}$ is a splitting
tree of $G^{*}$.

If both $C_{i_4}\cup v'v_{i_4}$ and $C_{2}\cup v_{2}v''\cup
v''v_{i_3}\cup C_{i_3}$ are connected component of $G^{*}$ which
contain an odd number of edges, then we will discuss it in the
following two subcases.

\texttt{\textit{Subcase 3.1.2-2a:}} \ In graph $G^{*}$, $C_{2}$ is a
connected component with an even number of edges, and $C_{i_3}$ is
one which contains an odd number of edges.

Noticing that both $C_{2}$ and  $C_{i_3}\cup v_{i_3}v''\cup
v''v'\cup v'v_{i_4}\cup C_{i_4}$ are connected component of $G^{*}$
which contain an even number of edges, we can easily get that $T^{*}
= (G^{*} \cap \mathbb{T}) \cup \{v'v_{1}, v_2v''\}$ is a splitting
tree of $G^{*}$, which is depicted in Fig.17.

\texttt{\textit{Subcase 3.1.2-2b:}} \ In graph $G^{*}$, $C_{2}$ is a
connected component with an odd number of edges, and $C_{i_3}$ is
one which contains an even number of edges.

Because there is exactly one $u,\omega$-$path$ in $\mathbb{T}$ for
any two vertices $u$ and $\omega$ in $G$, and none of \{$vv_2, vv_3,
vv_4$\} is an edge of $\mathbb{T}$, there must be exactly one
$v,v_3$-$path$ and exactly one $v,v_4$-$path$ in $\mathbb{T}$, and
the $v,v_3$-$path$, $v,v_4$-$path$ in $\mathbb{T}$ must be of the
form as $vv_1\dots v_{3}$ and  $vv_1\dots v_4$ respectively.
Noticing that, in the graph $G^{*}$, the connected components
$C_{i_3}$ and $C_{2}\cup v_{2}v''\cup v''v'\cup v'v_{i_4}\cup
C_{i_4}$ both contain an even number of edges, we can easily get
that $T^{*} = (G^{*} \cap \mathbb{T}) \cup \{v'v_{1}, v''v_{i_3}\}$
is a splitting tree of $G^{*}$, which is depicted in Fig.18.

\medskip

\textsf{Subcase 3.2:} \ One of \{$vv_3, vv_4$\} is the edge in
$\mathbb{T}$.

\medskip

Without loss of generality, let $vv_4$ be the edge in $\mathbb{T}$,
which is depicted in Fig.20. In addition, throughout Subcase 3.2,
let $vv_3$ and $vv_4$ be replaced by $v''v_{3}$ and $v'v_{4}$
respectively after the vertex splitting on $v$ in $G$; and the edge
set \{$vv_1, vv_2$\} be replaced by \{$v'v_{i_1}, v''v_{i_2}$\},
where \{$v_{i_1}, v_{i_2}$\}=\{$v_{1}, v_{2}$\} and \{$C_{i_1},
C_{i_2}$\}=\{$C_{1}, C_{2}$\}, which is depicted in Fig.19.
According to the edge $v_1v_2$ of $G$ is in the splitting tree
$\mathbb{T}$ or not, it will be discussed in the following two
subcases.

\begin{footnotesize}

\setlength{\unitlength}{0.6mm}
\begin{center}
\begin{picture}(70,36)

\put(-100,-3){\begin{picture}(30,40)

\put(16,28){\circle*{1.5}}

  \put(17,24){$v_{i_1}$}

   \put(6.3,30.3){$C_{i_1}$}

\multiput(16,28)(-0.8,1.8){5}{\circle*{0.5}} 

\put(16,12){\circle*{1.5}}

  \put(17,14){$v_{i_2}$}

  \put(10.5,5.5){$C_{i_2}$}

\multiput(16,12)(-2,-0.5){5}{\circle*{0.5}}  

\put(30,28){\circle*{1.5}}

  \put(28.5,30){$v'$}

\put(30,12){\circle*{1.5}}

  \put(29,6.5){$v''$}

\put(30,28){\line(-1,0){14}}  

\put(30,12){\line(-1,0){14}}   

\multiput(30,12)(0,2){9}{\circle*{0.5}}   

\multiput(30,28)(2,0){8}{\circle*{0.5}}  

\put(30,12){\line(1,0){14}}   

\put(44,28){\circle*{1.5}}

   \put(38.5,24){$v_{4}$}

\put(46,37){\circle*{1.5}}

   \put(47.5,37){$v_{l}$}

    \put(48,28){$C_{4}$}

\multiput(44,28)(0.6,2.2){4}{\circle*{0.5}}  

\multiput(46,37)(-1.8,0.6){5}{\circle*{0.5}}  

\put(44,12){\circle*{1.5}}

  \put(38,14){$v_{3}$}

\put(54,11){\circle*{1.5}}

  \put(56,9){$v_{p}$}

   \put(45.5,5){$C_{3}$}

\multiput(44,12)(1.8,-0.1){6}{\circle*{0.5}} 

\multiput(54,11)(0.6,1.8){5}{\circle*{0.5}}    

\put(8.5,35.5){\oval(21,20)[br]}

\put(8.5,4.5){\oval(21,20)[tr]}

\put(51.5,35.5){\oval(21,20)[bl]}

\put(51.5,4.5){\oval(21,20)[tl]}

\put(16,12){\line(0,1){16}} 

\put(19.5,-4){\scriptsize \bf {Fig.19}: $G^{*}$}

\end{picture}
}

\put(-30,-3){\begin{picture}(27,40)

\put(16,28){\circle*{1.5}}

  \put(19,27.5){$v_{1}$}

   \put(8.5,31){$C_{1}$}

\multiput(16,28)(-2,0.2){6}{\circle*{0.5}} 

\put(16,12){\circle*{1.5}}

  \put(19,10){$v_{2}$}

    \put(7,8){$C_{2}$}

\multiput(16,12)(-0.8,-1.8){6}{\circle*{0.5}} 

\put(30,20){\circle*{1.5}}

  \put(29,22){$v$}

\put(44,28){\circle*{1.5}}

   \put(37,28.5){$v_{4}$}

\put(46,37){\circle*{1.5}}

   \put(47.5,37){$v_{l}$}

    \put(48,28){$C_4$}

\multiput(44,28)(0.6,2.2){4}{\circle*{0.5}}  

\multiput(46,37)(-1.8,0.6){5}{\circle*{0.5}}  

\put(44,12){\circle*{1.5}}

  \put(37,10.5){$v_{3}$}

\put(54,11){\circle*{1.5}}

  \put(56,9){$v_{p}$}

   \put(46,4){$C_{3}$}

\multiput(44,12)(1.8,-0.1){6}{\circle*{0.5}} 

\multiput(54,11)(0.6,1.8){5}{\circle*{0.5}}    

\put(8.5,35.5){\oval(21,20)[br]}

\put(8.5,4.5){\oval(21,20)[tr]}

\put(51.5,35.5){\oval(21,20)[bl]}

\put(51.5,4.5){\oval(21,20)[tl]}

\put(30,20){\line(-5,3){14}}     

\put(30,20){\line(-5,-3){14}}   

\multiput(30,20)(1.7,1){8}{\circle*{0.5}} 

\put(30,20){\line(5,-3){14}}  

\put(16,12){\line(0,1){16}}    

\put(1,20){\vector(-1,0){9}}

\put(19.5,-4){\scriptsize \bf {Fig.20}: $G$}

\end{picture}
}
\put(35,-3){\begin{picture}(30,40)

\put(16,28){\circle*{1.5}}

  \put(17,24){$v_{i_1}$}

   \put(6.3,30.3){$C_{i_1}$}

\multiput(16,28)(-0.8,1.8){5}{\circle*{0.5}} 

\put(16,12){\circle*{1.5}}

  \put(17,14){$v_{i_2}$}

  \put(10.5,5.5){$C_{i_2}$}

\multiput(16,12)(-2,-0.5){5}{\circle*{0.5}}  

\put(30,28){\circle*{1.5}}

  \put(28.5,30){$v'$}

\put(30,12){\circle*{1.5}}

  \put(29,6.5){$v''$}

\put(30,28){\line(-1,0){14}}   

\put(30,12){\line(-1,0){14}}   

\put(30,12){\line(0,1){16}}    

\multiput(30,28)(2,0){8}{\circle*{0.5}} 

\multiput(30,12)(2,0){8}{\circle*{0.5}}   

\put(44,28){\circle*{1.5}}

   \put(38.5,24){$v_{4}$}

    \put(48,28){$C_{4}$}

\multiput(44,28)(0.9,2){5}{\circle*{0.5}}  

\put(44,12){\circle*{1.5}}

  \put(38.2,14){$v_{3}$}

   \put(45.5,4){$C_{3}$}

\multiput(44,12)(1.8,-0.5){6}{\circle*{0.5}} 

\put(8.5,35.5){\oval(21,20)[br]}

\put(8.5,4.5){\oval(21,20)[tr]}

\put(51.5,35.5){\oval(21,20)[bl]}

\put(51.5,4.5){\oval(21,20)[tl]}

\multiput(16,12)(0,2){9}{\circle*{0.5}} 

\put(19.5,-4){\scriptsize \bf {Fig.21}: $G^{*}$}

\end{picture}
}

\put(100.3,-3){\begin{picture}(30,40)

\put(16,28){\circle*{1.5}}

  \put(17,24){$v_{i_1}$}

   \put(6.3,30.3){$C_{i_1}$}

\multiput(16,28)(-0.8,1.8){5}{\circle*{0.5}} 

\put(16,12){\circle*{1.5}}

  \put(17,14){$v_{i_2}$}

  \put(10.5,5.5){$C_{i_2}$}

\multiput(16,12)(-2,-0.5){5}{\circle*{0.5}}  

\put(30,28){\circle*{1.5}}

  \put(28.5,30){$v'$}

\put(30,12){\circle*{1.5}}

  \put(29,6.5){$v''$}

\put(30,28){\line(-1,0){14}}   

\multiput(30,12)(-2,0){8}{\circle*{0.5}}   

\put(30,12){\line(0,1){16}}    

\multiput(30,28)(2,0){8}{\circle*{0.5}} 

\put(30,12){\line(1,0){14}}   

\put(44,28){\circle*{1.5}}

   \put(38.5,24){$v_{4}$}

    \put(48,28){$C_{4}$}

\multiput(44,28)(0.9,2){5}{\circle*{0.5}}  

\put(44,12){\circle*{1.5}}

  \put(38.2,14){$v_{3}$}

   \put(45.5,4){$C_{3}$}

\multiput(44,12)(1.8,-0.5){6}{\circle*{0.5}} 

\put(8.5,35.5){\oval(21,20)[br]}

\put(8.5,4.5){\oval(21,20)[tr]}

\put(51.5,35.5){\oval(21,20)[bl]}

\put(51.5,4.5){\oval(21,20)[tl]}

\multiput(16,12)(0,2){9}{\circle*{0.5}} 

\put(19.5,-4){\scriptsize \bf {Fig.22}: $G^{*}$}

\end{picture}
}

\end{picture}
\end{center}
\end{footnotesize}

\medskip

\textsf{\textit{Subcase 3.2.1:}} \ In graph $G$, $v_1v_2$ is not an
edge of $\mathbb{T}$.

\medskip

In this case, it is obvious that $T^{*} = (G^{*} \cap \mathbb{T})
\cup \{v'v_{4}, v'v''\}$ is a splitting tree of $G^{*}$, which is
depicted in Fig.19. So, in Subcase 3.2.1, $G^{*}$ is upper
embeddable.

\medskip

\textsf{\textit{Subcase 3.2.2:}} \ In graph $G$, $v_1v_2$ is an edge
of $\mathbb{T}$.

\medskip

In this case, if $C_{i_1}$ in $G^{*}$ is the same connected
component with at least one of \{$C_{i_2}, C_{3}, C_{4}$\}, then
$T^{*} = (G^{*} \cap \mathbb{T}) \cup \{v'v_{4}, v'v''\}$ is a
splitting tree of $G^{*}$. If any pair of components, which is
selected from \{$C_{i_1}, C_{i_2}, C_{3}, C_{4}$\}, is not the same
connected component of $G^{*}$, then it will be discussed in the
following two subcases.

\texttt{\textit{Subcase 3.2.2-1:}} \ In graph $G^{*}$, $C_{i_1}\cup
v'v_{i_1}$ contains an even number of edges.

In this case, it is obvious that $T^{*} = (G^{*} \cap \mathbb{T})
\cup \{v'v_{4}, v'v''\}$ is a splitting tree of $G^{*}$,  So, in
Subcase 3.2.2-1, $G^{*}$ is upper embeddable.

\texttt{\textit{Subcase 3.2.2-2:}} \ In graph $G^{*}$, $C_{i_1}\cup
v'v_{i_1}$ contains an odd number of edges, and $C_{i_2}\cup
v_{i_2}v''\cup v''v_{3}\cup C_{3}$ contains an even number of edges.

In this case, the connected component $C_{1}\cup v_{1}v\cup
C_{2}\cup v_{2}v\cup vv_{3}\cup C_{3}$, which contains an odd number
of edges in $G$, is replaced by $C_{i_1}\cup v'v_{i_1}$  and
$C_{i_2}\cup v_{i_2}v''\cup v''v_{3}\cup C_{3}$ after the vertex
splitting on $v$ in $G$. Let $T^{*} = (G^{*} \cap \mathbb{T}) \cup
\{v'v_{4}, v'v''\}$. Because $C_{i_1}\cup v'v_{i_1}$ contains an odd
number of edges, and $C_{i_2}\cup v_{i_2}v''\cup v''v_{3}\cup C_{3}$
contains an even number of edges, it is obvious that $\xi(
G^{*},T^{*})$=$\xi(G, \mathbb{T}))\leqslant 1$. So, $T^{*}$ is a
splitting tree of $G^{*}$.

\texttt{\textit{Subcase 3.2.2-3:}} \ In graph $G^{*}$, both
$C_{i_1}\cup v'v_{i_1}$ and $C_{i_2}\cup v_{i_2}v''\cup v''v_{3}\cup
C_{3}$ contain an odd number of edges.

In this case, according to the parity of the number of the edges in
$C_{i_2}$ and $C_{3}$ respectively, it will be discussed in the
following two subcases.

\texttt{\textit{Subcase 3.2.2-3a:}} \ In graph $G^{*}$, $C_{i_2}$
contains an odd number of edges, and $C_{3}$ contains an even number
of edges.

Because there is exactly one $u,\omega$-$path$ in $\mathbb{T}$ for
any two vertices $u$ and $\omega$ in $G$, and none of \{$vv_1, vv_2,
vv_3$\} is an edge of $\mathbb{T}$, there must be exactly one
$v,v_3$-$path$ in $\mathbb{T}$, and the $v,v_3$-$path$ in
$\mathbb{T}$ must be of the form as $vv_4\dots v_{3}$. Noticing
that, in the graph $G^{*}$, the connected components $C_{3}$ and
$C_{i_2}\cup v_{i_2}v''\cup v''v'\cup v'v_{i_1}\cup C_{i_1}$ both
contain an even number of edges, we can easily get that $T^{*} =
(G^{*} \cap \mathbb{T}) \cup \{v'v_{4}, v''v_{3}\}$ is a splitting
tree of $G^{*}$, which is depicted in Fig.21.

\texttt{\textit{Subcase 3.2.2-3b:}} \ In graph $G^{*}$, $C_{i_2}$
contains an even number of edges, and $C_{3}$ contains an odd number
of edges.

In this case, noticing that in the graph $G^{*}$ the connected
components $C_{i_2}$ and $C_{3}\cup v_{3}v''\cup v''v'\cup
v'v_{i_1}\cup C_{i_1}$ both contain an even number of edges, we can
easily get that $T^{*} = (G^{*} \cap \mathbb{T}) \cup \{v'v_{4},
v''v_{i_2}\}$ is a splitting tree of $G^{*}$, which is depicted in
Fig.22.

\medskip

From Case 1, Case 2, and Case 3, the Lemma 2.2 is obtained.
$\hspace*{\fill} \Box$

\bigskip

{\bf Theorem  2.1}   \ \ Let $G$ be a graph with minimum degree at
least 3, $v$ be a vertex of $G$ with deg$_{G}(v)$ $\geqslant$ 4,
$G^{*}$ be the graph obtained from $G$ by splitting $v$ into two
adjacent vertices $v'$ and $v''$, furthermore, the $v$-$local$
subgraph $G_{loc}(v)$ be connected. Then the graph $G$ is upper
embeddable if and only if $G^{*}$ is upper embeddable.

\medskip

{\bf Proof } \ \ ($\Longleftarrow$) \ \ \  Let $\textit{E}^{*}$ be
an embedding of $G^{*}$ in the orientable surfaces $S_{g}$ of genus
$g$. Then we can get an embedding $\textit{E}$ of $G$ in the surface
$S_{g}$ by contracting the $splitting$-$edge$ $v'v''$ in
$\textit{E}^{*}$. So
$\lfloor\frac{\beta(G)}{2}\rfloor$=$\lfloor\frac{\beta(G^{*})}{2}\rfloor$=$\gamma_{M}(G^{*})\leqslant
\gamma_{M}(G)$. On the other hand, $\gamma_{M}(G)\leqslant
\lfloor\frac{\beta(G)}{2}\rfloor$. Therefore,
$\gamma_{M}(G)=\lfloor\frac{\beta(G)}{2}\rfloor$, $i.e.$, the graph
$G$ is upper embeddable.

($\Longrightarrow $) \ \ \ Let $v_{1}$, $v_{2}$, $\dots$,  $v_{n}
(n\geqslant4)$ be all the vertices adjacent to $v$ in $G$, $v'$and
$v''$ be the replacement of $v$ after the vertex splitting on $v$ in
$G$, and the edge subset $\{vv_{i}|i=1,2,\dots n\}$ of $E(G)$ is
replaced by the subset $\{v^{*}v_{i}|v^{*}$ may be $v'$ or $v''$,
$i=1,2,\dots n\}$ of $E(G^{*})$. It can be obtained from Lemma 2.1
that there exists a splitting tree $\mathbb{T}$ of $G$ such that all
of \{$vv_1, vv_2, \dots, vv_{n}$\} are edges of $\mathbb{T}$. Let
$T^{*}=\{G^{*}\cap \mathbb{T}\}\cup v'v'' \cup \{v^{*}v_{i}|$$v^{*}$
may be $v'$ or $v''$, $i=1,2,\dots n\}$. Obviously, $T^{*}$ is a
spanning tree of $G^{*}$, and $\xi(G^{*},T^{*})$ =
$\xi(G,\mathbb{T}) = \xi(G) \leqslant 1$. So $T^{*}$ is a splitting
tree of $G^{*}$, and $G^{*}$ is upper embeddable. $\hspace*{\fill}
\Box$

\medskip

Especially, for a vertex $v$ of $G$ with deg$_{G}(v)$=4, we have the
following theorem.

\bigskip

{\bf Theorem  2.2}   \ \ Let $G$ be a graph with minimum degree at
least 3, $v$ be a vertex of $G$ with deg$_{G}(v)$=4, $G^{*}$ be the
graph obtained from $G$ by splitting $v$ into two adjacent vertices
$v'$ and $v''$, where the $splitting$-$edge$ $v'v''$ is not a
cut-edge of the $v$-$splitting$ subgraph $G^{*}_{spl}(v)$. Then the
graph $G$ is upper embeddable if and only if  $G^{*}$ is upper
embeddable.

\medskip
{\bf Proof } \ \ ($\Longleftarrow $) \ \ \ It is the same with that
of the Theorem 2.1.

($\Longrightarrow $) \ \ \ It is an obvious result of the Lemma 2.2.
$\hspace*{\fill} \Box$

 \bigskip

\noindent {\bf 3. Weak minor and upper embeddability}

\bigskip
In this section, we will provide a method to construct a
weak-minor-closed family of upper embeddable graphs from the bouquet
of circles $B_{n}$; in addition, we provide a corollary which
extends a result obtained by L. Nebesk\'{y} \cite{neb2}.

Let $v$ be a vertex of the graph $G$ with deg$_{G}(v)\geqslant 4$,
$G^{*}$ be the graph obtained from $G$ by splitting $v$ into two
adjacent vertices $v'$ and $v''$, then  $v$ is referred to as a
$flexible$-$vertex$ of $G$ if it satisfies one of the following two
conditions: (I) \ If $v$ is a vertex of the graph $G$ with
deg$_{G}(v)\geqslant4$, then the $v$-$local$ subgraph $G_{loc}(v)$
is connected (and the vertex splitting operation on this kind of
vertices is referred to as $type$-$I$ $vertex$ $splitting$); (II) \
If $v$ is a vertex of the graph $G$ with deg$_{G}(v)$=4, then the
$splitting$-$edge$ $v'v''$ is not a cut-edge of the $v$-$splitting$
subgraph $G^{*}_{spl}(v)$ (this kind of vertex splitting operation
is referred to as $type$-$II$ $vertex$ $splitting$).

According to Theorem 2.1 and Theorem 2.2, we can get, from the
bouquet of circles $B_{n}$, a weak-minor-closed family of upper
embeddable graphs through a sequence of vertex splitting operations
on the $flexible$-$vertices$.

\medskip

A graph $G$ is called locally connected if for every vertex $v$ of
$G$ the $v$-$local$ subgraph $G_{loc}(v)$ is connected.  In 1981, L.
Nebesk\'{y} \cite{neb2} obtained that every connected, locally
connected graph is upper embeddable. The following corollary extends
this result.

\bigskip

{\bf Corollary }   \ \ A graph, which is obtained from a connected,
locally connected graph through a sequence of type-I or type-II
vertex splitting operations on it, is upper embeddable.

\medskip

{\bf Proof } \ \ According to the result obtained  by L. Nebesk\'{y}
\cite{neb2} we can get that every connected, locally connected graph
is upper embeddable. Combining with Theorem 2.1 and Theorem 2.2 we
can get the Corollary.   $\hspace*{\fill} \Box$

\bigskip

\noindent {\bf 4. Conclusions}

\bigskip

{\bf Remark 1} \ \ Let $G$ be an upper embeddable graph with minimum
degree at least 3, $v$ be a vertex of $G$ with deg$_{G}(v) \geqslant
5$, $G^{*}$ be the graph obtained from $G$ by splitting $v$ into two
adjacent vertices $v'$ and $v''$. Then the condition that the
$splitting$-$edge$ $v'v''$ is not a cut-edge of the $v$-$splitting$
subgraph $G^{*}_{spl}(v)$ can not guarantee the upper embeddability
of $G^{*}$. For example, the graph $G^{*}$ in Fig.24 is a graph
obtained from the upper embeddable graph $G$ in Fig.23 through
vertex splitting on $v$ in $G$, and the $splitting$-$edge$ $v'v''$
is not a cut-edge of the $v$-$splitting$ subgraph $G^{*}_{spl}(v)$.
But,  $G^{*}$ is not upper embeddable.

\medskip

\begin{footnotesize}

\setlength{\unitlength}{0.6mm}
\begin{center}
\begin{picture}(100,18)


\put(-20,-13){\begin{picture}(10,10)

\put(0.5,32.6){\circle{7}}

\put(3,30){\circle*{1.5}}

\put(3,30){\line(0,-1){15}}

\put(3,30){\line(1,0){20}}

\put(24,30){\circle*{1.5}}

\put(32,36){\circle*{1.5}}

\put(32,24){\circle*{1.5}}

\put(32,24){\line(0,1){12}}

\put(24,30){\line(4,3){8}}

\put(24,30){\line(4,-3){8}}

\put(3,15){\circle*{1.5}}

\put(24,15){\circle*{1.5}}

\put(3,15){\line(1,0){20}}

\put(24,15){\line(0,1){15}}

\qbezier(3,15)(17,22)(24,15)

\qbezier(3,15)(17,8)(24,15)

\qbezier(3,15)(8,27)(24,30)

\put(21,32){$v$}

\put(4,3){\scriptsize{\bf Fig.23:} \ $G$}
\end{picture}}


\put(75,-13){\begin{picture}(10,10)

\put(0.5,32.6){\circle{7}}

\put(3,30){\circle*{1.5}}

\put(3,30){\line(0,-1){15}}

\put(3,30){\line(1,0){20}}

\put(24,30){\circle*{1.5}}

\put(24,23){\circle*{1.5}}

\put(32,36){\circle*{1.5}}

\put(32,26){\circle*{1.5}}

\put(32,26){\line(0,1){10}}

\put(24,30){\line(4,3){8}}

\put(24,30){\line(2,-1){8}}

\put(3,15){\circle*{1.5}}

\put(24,15){\circle*{1.5}}

\put(3,15){\line(1,0){20}}

\put(24,15){\line(0,1){15}}

\qbezier(3,15)(17,22)(24,15)

\qbezier(3,15)(17,8)(24,15)

\qbezier(3,15)(8,27)(24,23)

\put(21,32){$v'$}

\put(25.8,20){$v''$}

\put(4,3){\scriptsize{\bf Fig.24:} \ $G^{*}$}

\put(-40,20){\vector(1,0){20}}

\end{picture}}
\end{picture}
\end{center}

\end{footnotesize}

\bigskip

{\bf Remark 2} \ \ Let $v_1v_2$ be an edge of the graph $G$. The
$edge$-$global$ $subgraph$ of $v_1v_2$, which is denoted by
$G_{glo}(v_1v_2)$, is the subgraph of $G$ that is induced by the
vertices of $v_1$, $v_2$ and all the neighbors of them. The
$edge$-$local$ $subgraph$ of $v_1v_2$, which is denoted by
$G_{loc}(v_1v_2)$, is the subgraph of $G$ that is induced by all the
neighbors of the vertex $v_1$ and $v_2$. A $flexible$-$edge$ of
graph $G$ is such an edge $v_1v_2$ of $G$ which satisfies one of the
following two conditions: (I) $v_1v_2$ is not a cut-edge of the
$edge$-$global$ $subgraph$ of $v_1v_2$, and the adjacent vertices
$v_1$, $v_2$ are replaced by a vertex $v$ of degree 4 after
contracting the edge $v_1v_2$; (II) The $edge$-$local$ $subgraph$
$G_{loc}(v_1v_2)$ of $v_1v_2$ is connected, and the adjacent
vertices $v_1$, $v_2$ are replaced by a vertex $v$ with degree no
less than 4 after contracting the edge $v_1v_2$. A
$flexible$-$weak$-$minor$ of the graph $G$ is a graph obtained from
$G$ through a sequence of edge-contraction operations on the
$flexible$-$edges$.

From Theorem 2.1 and Theorem 2.2 we can get that a graph $G$ is
upper embeddable if and only if its $flexible$-$weak$-$minor$ is
upper embeddable. So the determining of the upper embeddability of
$G$ can be replaced by determining the upper embeddability of its
$flexible$-$weak$-$minor$. Furthermore, the algorithm complexity of
determining the upper embeddability of  $G$ may be reduced much by
this way, because the order of the $flexible$-$weak$-$minor$ of $G$
is less than the order of $G$.

 
\end{document}